\magnification=1200

\def\Z{{\bf Z}}
\def\R{{\bf R}}
\def\C{{\bf C}}
\def\Q{{\bf Q}}
\def\O{{\cal O}}
\def\F{{\cal F}}
\def\G{{\Gamma}}
\def\g{{ \bf g}}
\def\k{{\bf k}}
\def\ra{{\longrightarrow}}
\def\p{{\partial}}
\def\db{{\overline \p}}

\centerline{\bf ROZANSKY-WITTEN INVARIANTS VIA FORMAL GEOMETRY}\par
\medskip
\medskip
\centerline {\bf Maxim Kontsevich}

\vskip 15truemm

\noindent { 0. INTRODUCTION}
\vskip 0.5truecm

Recently  L.~Rozansky and E.~Witten (see [RW]) proposed a topological quantum 
field theory
 depending on a compact oriented 3-dimensional manifold $M$ and on a compact 
 hyperk\"ahler manifold 
 $X$.
 In the case 
 $H^1(M,\Q)=0$ 
 (i.e. when $M$ is a rational homology 3-sphere), the 
 partition function $Z(M,X)\in\C$ 
 of this field theory can be calculated using finitely many
  terms of the perturbation theory. 
  As a function on $M$ this is so called invariant of finite type,
  of  order $2n$ where $4n$ is the dimension of $X$. 
  
  More generally, for every connected finite $3$-valent graph $\G$ with $2n$ 
  vertices, endowed with a cyclic order in the star of each vertex, Rozansky
   and Witten
   associated a function $X\mapsto Z_{\G}(X)$
   on the space of isometry classes of hyperk\"ahler 
   manifolds of dimension $4n$. This function is given by an integral over $X$
   of certain invariant polynomial in coefficients of the curvature tensor
    of the hyperk\"ahler metric on $X$. 
    
  Here we propose   a simple construction of RW invariants.
   It consists of two  steps:
   
   1) with every $3$-valent graph endowed with orientations at vertices, 
   or with 
   every oriented rational homology
   $3$-sphere $M$, we associate a cohomology class of the Lie algebra of 
    formal Hamiltonian vector fields in an arbitrary finite-dimensional 
     symplectic vector space. This cohomology class
     is stable. Stable cohomology groups under the question are called
      Graph Cohomology groups, because they can be calculated via certain 
      complex
       constructed from finite graphs. Universal finite type
       invariants of links and of rational homology $3$-spheres take
        values in certain subspace of the Graph Cohomology.
       Morally, this construction is a universal perturbative 
       quantum Chern-Simons theory.
       
       This construction is known already for a while, see [Ko1] for
        the general overview, [BN] for the discussion of $3$-valent graphs
         and finite type invariants of links, 
         and [LT] for the construction of invariants of homology $3$-spheres.
       
   2) it is known from early 70-ies that
    cohomology groups of Lie algebras of formal
    vector fields give
    characteristic classes of foliations (see [BR], [BH]). 
    In the case of Hamiltonian vector
     fields
     we get characteristic classes of symplectic foliations, i.e. of foliations
      endowed with a symplectic structure in the transversal direction.
      Formally, any complex manifold can be considered as a foliation
       in anti-holomorphic direction. Analogously, a holomorphic symplectic
        manifold $X$ gives a formal complex-valued symplectic foliation.
     Applying a small modification of the original construction of 
     characteristic
     classes we get a homomorphism from the Lie algebra cohomology
      to the the cohomology of coherent sheaves $H^{\bullet}(X,\O_X)$.
       This construction could have been invented
        25 years ago.

       RW invariants can be read from these characteristic classes.
        Moreover, it follows directly from our description that
         one doesn't need the hyperk\"ahler metric on $X$ in the construction.
         One can formulate it in purely holomorphic terms
           (or in algebro-geometric terms if $X$ is algebraic).
         An easy argument shows that numbers $Z(M,X)$ or $Z_{\G}(X)$
      are deformation invariants of hyperk\"ahler manifolds $X$.
      
      As a by-product, we obtain a construction of finite type invariants
       of $3$-manifolds  based on symplectic foliations instead
        of hyperk\"ahler manifolds.

      This paper is an extended version of my letters to V.~Ginzburg and
      to E.~Witten (January 1997). Recently M.~Kapranov, stimulated by these
       letters, found a different  approach to RW invariants.  
      He noticed  that cohomology classes in $H^{\bullet}(X,\O_X)$
       associated with 
      all $3$-valent graphs
       can be written down in terms of just one class, so called Atiyah class.
        His construction is  shorter than mine, but basically  is
         the same. 
        M.~Kapranov wrote a  beautiful and detailed exposition (see [Ka])
         with many 
        interesting deviations from the  main theme.
        Still, I think that it is reasonable to give an account of 
      the original geometric approach. Strictly speaking, my
      present paper contains no really new ideas. Nevertheless,
      I hope that it could help to clarify 
      the picture.

      I am grateful to M.~Kapranov and to A.~Schwarz for useful communications,
       and to E.~Witten for showing me a preliminary version of the
        paper [RW].

\vskip 0.5truecm    
 \noindent { 0.1. {\it Notations}}   
\vskip 0.5truecm

  Let $\g$ be a Lie algebra over $\R$, $\k\subset \g$ be a finite-dimensional
  Lie subalgebra of $\g$ 
   and $K$ be a Lie group (not necessarily connected)
   with the Lie algebra $\k$. We assume also that an action of $K$ by 
   automorphisms of $\g$ is given, such that the induced action of $\k=Lie(K)$
   is the adjoint action. Let $V$ be a $(\g,K)$-module.
    Relative cohomology group $H^{\bullet}(\g,K;V)$ is defined as the cohomology group
    of the complex of $K$-invariant skew-symmetric
    polylinear maps from $\g$ to $V$ vanishing if one of arguments belongs to 
    $\k$:
    $$C^i(\g,K;V):=\left(Hom(\wedge^i(\g/\k), V\right)^K,\,\,\,i\ge 0\,.$$
   
    This complex is a subcomplex of the standard cochain complex of $\g$ with
     coefficients in $V$. The differential in $C^{\bullet}(\g,K;V)$ is induced
      from the standard differential in $C^{\bullet}(\g;V)$.
     
     For the case of trivial coefficients the cochain
     complex $C^{\bullet}(\g;\R)$ can be interpreted as the complex of
      (left) $G$-invariant differential forms on $G$, where $G$ is any Lie 
      group with the Lie algebra equal to $\g$. Analogously, the 
      relative cochain 
      complex $C^{\bullet}(\g,K;\R)$ can be identified with the
       complex of $G$-invariant differential forms on $G/K$.

   For $n\ge 0$ we denote by $Ham_{2n}$ the Lie algebra of formal Hamiltonian
    vector fields in the standard symplectic vector space $\R^{2n}$.
     This Lie algebra is endowed with the topology of the inverse limit
      of finite-dimensional vector spaces.
    Elements of $Ham_{2n}$ are in one-to-one correspondence with
    formal  Hamiltonians modulo constants:
    $$Ham_{2n}\simeq \R[[p_1,\dots,p_n,q_1,\dots,q_n]]/\R\,\,\,.$$
    
   We denote by $Ham^0_{2n}$ the subalgebra of $Ham_{2n}$ consisting
    of formal vector fields vanishing at zero:
      $$Ham_{2n}/Ham^0_{2n}\simeq \R^{2n}\,\,\,.$$
      
   Lie algebra $Ham^0_{2n}$ contains subalgebra $sp(2n,\R)$ consisting 
   of linear
   Hamiltonian vector fields. The Lie group $Sp(2n,\R)$ acts on $Ham^0_{2n}$.
    Thus, we can define cohomology groups with coefficients in the trivial 
    one-dimensional module
    $$H^i_{2n}:=H^i_{cont} (Ham^0_{2n},Sp(2n,\R);\R)\,\,\,.$$
    
   Here the subscript $cont$ means that we consider only
    continuous cochains, 
   i.e. polylinear  functionals depending  on finitely many terms
    in the Taylor expansions at zero (Gelfand-Fuks cohomology).

\vskip 0.5truecm    
 \noindent { 1. FROM GRAPHS TO COHOMOLOGY}    
\vskip 0.5truecm

Let $\G$ be a  $3$-valent graph with $2N$ vertices. We  associate with it 
 an element $I_\G$
  in $H^{2N}_{2n}$ for any $n$. 
  
  Lie algebra $Ham^0_{2n}$ is a semi-direct product of $sp(2n,\R)$
   and of the subalgebra $Ham^1_{2n}$  consisting
    of Hamiltonians $H\in \R[[p_i,q_i]]$ such that the Taylor series of $H$
     starts at terms of order at least $3$. Thus, the relative cochains
      $C^i_{cont}(Ham^0_{2n},Sp(2n,\R);\R)$ can be identified with 
      $Sp(2n,\R)$-invariant 
    cochains
       of $Ham^1_{2n}$ with trivial coefficients.
     
     The group $Sp(2n,\R)$ acts semi-simply on $Ham^1_{2n}$, and on its cochain
     complex.
      This  implies that we have an isomorphism
      $$H^i_{2n}\simeq \left(H^i_{cont}(Ham^1_{2n};\R)\right)^{Sp(2n,\R)}
      \,\,\,.$$
      
      The first cohomology group of $Ham^1_{2n}$ (i.e. the co-abelianization)
      is non-trivial, 
      it contains $Sym^3(\R^{2n})$. The corresponding $1$-cochain associates
       to a formal Hamiltonian $H$ its third Taylor coefficient.
        Using cup-products we  construct cohomology 
        classes in higher degrees:
        $$\wedge^{2N}(Sym^3(\R^{2n}))\ra H^{2N}_{cont}(Ham^1_{2n};\R)\,\,\,.$$
        
        The map from above is evidently $Sp(2n,\R)$-equivariant. Now,
         any $3$-valent graph $\G$ with $2N$ vertices (the number of vertices
         of $\Gamma$ is necessarily even)
        gives, up to a sign, an invariant tensor in 
        $\wedge^{2N}(Sym^3(\R^{2n}))$. We use graph as a scheme for 
        contracting indices. Applying isomorphisms as above we get
          relative cohomology class $I_\G$. 
         
         For general discussion of Graph Cohomology and its relation
          with cohomology of Lie algebras of formal vector fields
          we refer the reader to [Ko2].

 \vskip 0.5truecm    
 \noindent { 2. CHARACTERISTIC CLASSES OF FLAT BUNDLES }    
\vskip 0.5truecm

 Before going further  we  remind a general construction of characteristic 
 classes of flat bundles. In Section 3 we will apply it to foliations.
 
 Let $X$ be a smooth manifold, $G$ be a finite-dimensional Lie group with the 
 Lie algebra $\g$. Let $E\ra X$ be a principal $G$-bundle endowed with a flat
  connection $\nabla$. 
  
  Assume that $E$ is trivial as a topological $G$-bundle. Let us choose
   a smooth trivialization of $E$. Then the connection $\nabla$ is given by
    a $1$-form $A$ on $X$ with values in $\g$, satisfying the Maurer-Cartan
     equation 
     $$dA+{1\over 2} [A,A]=0\,\,\,.$$
     
    We can consider $A$ as a linear map from $\g^*=C^1(\g;\R)$ to $\Omega^1(X)$.
   Let us extend it to the map from the whole cochain complex of $\g$
   $$\bigoplus_i \wedge^i(\g^*)=\bigoplus_i C^i(\g;\R)$$
   to $\bigoplus_i \Omega^i(X)$ using cup-products on the cochain complex and 
   on
    differential forms. The Maurer-Cartan equation guarantees  that this map 
     is a morphism of complexes. Thus, we have a map of cohomology groups:
     $$H^i(\g;\R)\ra H^i(X,\R)\,\,\,.$$
     
    This map does not change if we choose another trivialization of $E$
     in the same homotopy class of trivializations. The proof is immediate
     because  in such a situation we have a flat connection in the 
     trivialized $G$-bundle over the product $X\times[0,1]$.
     
     Another way to describe the same construction is to use the natural
    $G$-invariant $\g$-valued $1$-form $A_E$ on $E$ satisfying the 
    Maurer-Cartan equation. It gives a homomorphism  $H^i(\g,\R)\ra H^i(E,\R)$.
     A trivialization of $E$ gives a section $s:X\ra E$ of $E$. Then we can 
     restrict cohomology classes form $E$ to $s(X)$.

     For a semisimple group $G$ and for primitive classes in $H^i(\g;\R)$
      we get odd-dimensional characteristic classes for flat connections
      in topologically trivial bundles, essentially
       Chern-Simons secondary characteristic classes.
        For connected nilpotent groups $G$ we always get characteristic classes
         because these groups are
          contractible and all $G$-bundles are topologically
        trivial.
    
    Suppose now that we fixed a  Lie subgroup $K$ of $G$ such that
     the bundle $E$
    is topologically equivalent to a bundle induced from a $K$-bundle. 
     In other words, the bundle with fiber $G/K$ associated with $E$, 
      has a
     continuous  section. 
   Analogously to the previous construction, we can define a map
   $$C^i(\g,K;\R)\ra \Omega^i(X)\,\,\,.$$

  The induced map on cohomology 
  $$H^i(\g,K;\R)\ra H^i(X,\R)$$
  depends only on the homotopy class of an identification of $E$
    with induced bundles (i.e. of the section of the associated $G/K$-bundle).
    
   If the inclusion $K\subset G$ is  homotopy equivalence then 
   there is unique homotopy class of identifications. 
    For example, flat connections in complex $n$-dimensional vector bundles
     have characteristic classes (Chern-Simons classes)
     
     $$cs_i\in H^{2i+1}(gl(n,\C),U(n);\R)\ra H^{2i+1}(X,\R),\,\,\,0\le i\le n-1
     \,\,\,.$$
     
     The next generalization consists in consideration of an 
     infinite-dimensional
      group  $G$, for example a projective limit of finite-dimensional groups.
    Also, one can consider ``Lie groups'' which are formal manifolds
    in some directions. The algebra of functions on such a ``group''
     is the algebra of
    formal power series in several variables with coefficients
    in usual smooth functions on a Lie group. In algebraic geometry one
     can consider for these purposes a mixture of pro- and ind- schemes.
      Any pair $(\g,K)$ (as in the definition of relative cohomology, see 
      subsection 0.1) 
       produces a partially formal Lie group.

      An important example is the  ``Lie group'' $FDiff(\R^n)$ 
       (formal diffeomorphisms of $\R^n$). Its Lie algebra is
         the Lie algebra $W_n$
      of formal vector fields in $\R^n$. The underlying
        topological space of $FDiff(\R^n)$ 
        is the group all automorphisms of $\R$-algebra
         $\R[[x_1,\dots, x_n]]$ (i.e. the group of
         formal diffeomorphisms of $\R^n$
          fixing $0$). Functions on this group are
        formal power series on $n$-dimensional space 
        $$\R^n\simeq \R\langle {\p\over\p x_1},\dots,{\p\over\p x_n}\rangle$$
     with coefficients in smooth functions on spaces of
      jets of a finite order of diffeomorphisms
      of $\R^n$ fixing $0$.

 \vskip 0.5truecm    
 \noindent {3. FOLIATIONS AND FLAT BUNDLES }    
\vskip 0.5truecm

      If $X$ is a smooth manifold and $\F$ is a foliation of codimension $n$
       on $X$ then we construct a  linear map
        $$H^i_{cont}(W_n,O(n,\R);\R)\ra H^i(X,\R)\,\,\,.$$
        
        Namely, we have a principal bundle over $X$ with the structure group
          equal to the group of formal diffeomorphisms of $\R^n$ fixing $0$.
         The fiber of this bundle at each point
         $x\in X$ consists of identifications of formal neighborhoods
        of $x$ in the space of leaves 
        $U/\F$, where  $U$ is sufficiently small neighborhood of $x$
         in $X$, with the formal neighborhood of $0$ in $\R^n$.
         This bundle carries a natural flat connection only along  $\F$.
          We can consider the associate principal $FDiff(\R^n)$-bundle. This 
          bundle carries a natural flat connection along all directions in $X$.
         
         As a topological group $FDiff(\R^n)$ is the same as  
         $Aut(\R[[x_1,\dots,x_n]])$
          (because formal coordinates do not change the topology).
           The group $Aut(\R[[x_1,\dots,x_n]])$ is homotopy equivalent to 
           $GL(n,\R)$, being a pro-nilpotent extension of it.
            The group $GL(n,\R)$ is homotopy equivalent to its subgroup
             $O(n,\R)$. Thus, we get characteristic classes as above.

       What we explained above
        is the standard construction of characteristic classes
        of foliations phrased in somewhat new terms.
         
 \vskip 0.5truecm    
 \noindent {4. APPLICATION TO HAMILTONIAN VECTOR FIELDS}    
\vskip 0.5truecm

Let us suppose now that the foliation $\F$ is endowed with the transversal
 symplectic structure and has codimension $2n$ instead of $n$ as above.
  To have such a structure is the same as to have a closed degenerate
   $2$-form $\omega$
   on $X$ of constant rank $2n$. The foliation $\F$ is given by the kernel
    of $\omega$.

 \vskip 0.5truecm    
 \noindent { 4.1. {\it Remark}}    
\vskip 0.5truecm 
   
    As a side remark, we want to notice that
    the natural source of degenerate 2-forms is the variational principle. 
    On the set of solutions of Euler-Lagrange equations
     one has naturally a closed $2$-form, which could be degenerate in some 
     cases. The standard point of view is opposite to this. People
      usually consider Poisson manifolds, i.e.  bivector fields satisfying the 
      Jacobi
      identity
      as a degeneration of symplectic
       geometry. In general, Poisson manifolds describe a limiting behavior of 
        quantum 
       mechanics.
       \vskip 0.5truecm
    
  For a foliation $\F$ one can introduce the de Rham complex along $\F$:
  $$\Omega^i_\F(X):=\Gamma(X,\wedge^i(T^*_\F))$$
  where $T_\F$ denote the tangent bundle to $\F$. This complex is a quotient
   complex of the de Rham complex of $X$. Cohomology $H^i_\F(X)$ of
    the complex $\Omega^{\bullet}_\F(X)$ is quite a wild
    object, non-computable in simple terms in general.  
    
    As in Section 3, we have a principal bundle over $X$ with the 
    structure group
     $$Aut(\R[[p_i,q_i]],\sum dp_i\wedge dq_i)\,\,\,.$$
      This bundle carries a natural
      flat connection along $\F$. The structure group is homotopy equivalent
     to the subgroup $Sp(2n,\R)$, and also to a smaller subgroup
      $U(n)$. Hence, we have natural maps
      $$H^i_{2n}=H^i_{cont}(Ham^0_{2n},Sp(2n,\R);\R)\ra
       H^i_{cont}(Ham^0_{2n}, U(n);\R)\ra H^i_\F(X)\,\,\,.$$

       Foliation $F$ has a natural transversal volume element $vol$
       represented by the 
       differential form $\omega^n/n!$ where $\omega$ is the closed $2$-form
        defining $\F$. Multiplication by this volume element gives a morphism 
        of complexes and of cohomology spaces:
        $$vol\wedge: \Omega^i_\F(X)\ra \Omega^{i+2n}(X),\,\,\,\,H^i_\F(X)\ra 
        H^{i+2n}(X,\R)\,\,\,.$$
        
       Thus, we get characteristic classes with values in the 
       de Rham cohomology
        of $X$ in degrees shifted by $2n$. Alternatively,
         there is an analogous map for Lie algebra cohomology:
         $$vol\wedge: H^i_{2n}\ra H^{i+2n}_{cont}
         (Ham_{2n},Sp(2n,\R);\R)\,\,\,.$$

 \vskip 0.5truecm    
 \noindent {5. GROUPS OF SYMPLECTOMORPHISMS}    
\vskip 0.5truecm      
       
       Let $(Y,\omega)$ be a compact symplectic $2n$-dimensional manifold.
     Let us denote by $Sympl(Y)$ (or $Sympl(Y,\omega)$) the
     group of symplectomorphisms of $Y$, and by $sympl(Y)$ the 
     Lie algebra of Hamiltonian vector fields on $Y$. 
      A version of the construction of characteristic classes of foliations
       gives the following homomorphism:
     $$H^i_{2n}=H^i_{cont}(Ham^0_{2n},Sp(2n,\R);\R)\ra
       H^i_{cont}(Ham^0_{2n}, U(n);\R)\ra H^i(Sympl(Y)^{\delta},\R)\,\,\,.$$
   Here the upper index $\delta$ means that we consider $Sympl(Y)$ as a discrete
      group. 
      
      This  homomorphism 
        takes values in the group of characteristic classes 
        of flat non-linear symplectic
        bundles.
      Let $E\ra B$ be a smooth bundle with a flat connection and a covariantly
       constant symplectic form on fibers. We assume that  fibers are isomorphic
       as 
        symplectic manifolds to $Y$. Such a structure is given by
         a homomorphism 
        (up to a conjugacy)
   $$\pi_1(B)\ra Sympl(Y)^{\delta}$$
   The total space $E$ carries a symplectic foliation. By constructions
    described above in Section 4 
     we get classes
    with values in $H^{i+2n}(E,\R)$. Fibers of 
    the bundle $E\ra B$ are compact and 
     naturally oriented. Thus, we can integrate  cohomology classes
      along fibers landing at the space  $H^i(B,\R)$.

      Analogously, one can construct homomorphisms for the Lie algebra
      cohomology
       $$H^i_{2n}=H^i_{cont}(Ham^0_{2n},Sp(2n,\R);\R)\ra
       H^i_{cont}(Ham^0_{2n}, U(n);\R)\ra H^i(sympl(Y),\R)\,\,\,.$$ 
      We refer the reader to [F].
        In general, it seems that there is a natural homomorphism from
      $H^i_{cont}(Ham^0_{2n}, U(n);\R)$  to the Van Est
       cohomology of $Sympl(Y)$ (the cohomology of a subcomplex of continuous
         cochains of the group $Sympl(Y)$).
        
  \vskip 0.5truecm    
 \noindent { 5.1. {\it Example of a flat symplectic bundle}}    
\vskip 0.5truecm       
         
         There is a natural series of finite-dimensional bundles
       with flat symplectic connections. The base is the moduli space of
        complex curves of genus $g\ge 2$, the fiber is the
         moduli space of irreducible
          unitary local systems (or, generally, flat connections
           with compact structure groups)
          on a surface of genus $g$. There is a standard symplectic
           form on the moduli space of flat connections which is defined
           purely topologically. Of course, sometimes such moduli spaces
            are non-compact, or singular after the compactification.
              Modulo these technical difficulties,
               this example gives a series of homomorphisms, labeled by
               compact Lie groups, from
              the Graph Cohomology to the cohomology of moduli spaces of 
              curves.
           In [Ko2] we constructed another homomorphism which maps
            the Graph Cohomology
            to the {\it homology \rm} groups of moduli spaces of curves.

         \vskip 0.5truecm    
 \noindent {6. COMPLEX GEOMETRY AND FOLIATIONS}    
\vskip 0.5truecm
         
         Let $X$ be a complex manifold of dimension $N$. 
         We denote by ${\widetilde X}$ the underlying
          smooth manifold of dimension $2N$. It is well-known that
           the almost-complex structure of $X$ can be considered as a 
           vector subbundle $T^{0,1}$ of the complexified tangent bundle 
           $T_{\widetilde X}\otimes \C$. The integrability of 
           almost complex structures is equivalent to the formal integrability
           of $T^{0,1}$. Thus, we get formally  a ``complex foliation'' on 
           ${\widetilde X}$. 
           
           There is still a better point of view.
            In order to describe it we introduce an auxiliary
            ``complex manifold'' ${\widetilde X}_\C$. 
             The underlying topological space of 
           this manifold ${\widetilde X}_\C$ is $\widetilde X$.
          The sheaf of functions on ${\widetilde X}_\C$ is the sheaf 
          of complex-valued smooth functions on ${\widetilde X}$ considered
          as an algebra over $\C$. We look at this sheaf as at a completion
         of the sheaf of holomorphic functions on complex manifold
          $X\times {\overline X}$ defined in small neighborhoods of the closed
         subset $X_{diag}:=\{(x,x)|x \in X\}$ of $X\times{\overline X}$.
          From this picture 
          it is clear that ${\widetilde X}_\C$ has formally the 
          structure of the product of two manifolds, and carries two
           transversal foliations. We would like to forget about one of them
         and leave the other. Thus, holomorphic functions on $X$ are functions
        constant along $\db$-foliation. Also, the de Rham complex along
         $\db$-foliation is nothing but the standard Dolbeault complex
          of $X$. Absence of higher cohomology groups for ``coherent sheaves''
           on ${\widetilde X}_\C$ can be viewed as a Stein property.
            In algebraic geometry we would call such spaces affine schemes.
         
         Let us return to symplectic geometry. If $X$ carries a holomorphic
         symplectic form $\omega$,  then on ${\widetilde X}_\C$ we have a 
         holomorphic
         symplectic 
         foliation. Now we can apply the same construction as in Section 4 
          and get a map
      $$H^i_{2n}\otimes \C:=H^i_{cont} (Ham^0_{2n}\otimes \C,Sp(2n,\C);\C)
      \ra H_\db^i ({\widetilde X})=H^i(X,\O)\,\,\,.$$
      
      It is almost evident from our description that this construction is
       complex-analytic, i.e. if $(X,\omega)$ holomorphically depends
        on parameters then corresponding classes also depend holomorphically.
      Moreover, the construction can be phrased in the language of
       algebraic geometry, see [Ka].

       A small extension of this construction involves cohomology with 
       non-trivial coefficients. For example, we have natural maps
      $$H^i_{cont} (Ham^0_{2n}\otimes \C,Sp(2n,\C);\wedge^j(\C^{2n}))
      \ra H^i(X,\Omega^j)\,\,\,.$$
      Corresponding graph complexes are associated with graphs with free legs.
       Part of these graph cohomology spaces corresponding to $3$-valent graphs
       appears as universal Vassiliev invariants of knots, see [BN].
      Definitely there are other non-trivial
 cohomology classes corresponding to graphs of higher valency,
 as follows from simple estimates of Euler characteristics of Graph Complexes.
       
       The construction of M.~Kapranov of characteristic classes
        associated with 3-valent graphs can be phrased as follows.
         The Atiyah class $\alpha_T$,
         introduced in [Ka], is the image in
         $H^1(X, Sym^3 T_X)$ of a natural class (see Section 1) 
          in 
          $$H^1_{cont}(Ham^0_{2n}\otimes \C,Sp(2n,\C);Sym^3(\C^{2n}))\,\,\,.$$
           Together with the symplectic form $\omega$, 
           which is an element of
            $H^0(X, \wedge^2(T^*_X))$ (or, of $H^0_{cont}
            (Ham^0_{2n}\otimes \C,Sp(2n,\C);\wedge^2(\C^{2n}))\,\,\,$), one can 
            construct characteristic classes contracting indices in the
      tensor product of copies of $\alpha_T$ and of $\omega$.
       
   \vskip 0.5truecm    
 \noindent {7. HYPERK\"AHLER MANIFOLDS}    
\vskip 0.5truecm    
       
       If $X$ is a compact hyperk\"ahler manifold then we have have 3 complex
       structures $I,J,K$ on $X$. Let us pick one of them, say  $I$.
        Complex manifold $X_I$ carries a holomorphic symplectic form $\omega_I$.
       We can construct numerical invariants
        multiplying characteristic classes 
        of $(X_I,\omega_I)$
         in $H^i(X_I,\Omega^j)$ by appropriate
         powers of the holomorphic symplectic form and of the cohomology
          class of the K\"ahler form, and then integrating over $X$.
           In [Ka] the reader can find  arguments showing that we get the same
           formulas as in the paper of Rozansky and Witten.
          
          For $3$-valent graphs the number which we get is invariant under
         deformations preserving the cohomology class of the K\"ahler form.
          The argument is is that 1) $Z_{\Gamma}(X,\omega)$
           it is a holomorphic function
           on the moduli space of complex symplectic manifolds
            with fixed polarization, depending only on the symplectic form
             modulo the multiplication by a constant scalar,
            2) by twistor construction one can produce a lot of rational
          in these moduli spaces, 3) 
          holomorphic functions on $\C P^1$ are constant.  
          
          The moduli space of complex structures on a compact
          hyperk\"ahler manifold is a locally symmetric Hermitean spaces
          of non-compact type. Thus, hyperk\"ahler manifolds should
           have degenerations to (possibly) simpler objects.
           Eventually, one expects that one can get a combinatorial
            objects at the limit, something like toric varieties. These
             combinatorial objects should produce weight systems for
              Vassiliev invariants.

 \vskip 0.5truecm    
 \noindent {8. SUPERSYMMETRIC FORMULATION}    
\vskip 0.5truecm

We have seen two sources of linear functionals on Graph Cohomology
 (and invariants of knots and rational homology $3$-spheres):
  symplectic foliations and complex symplectic manifolds.
   Another  (standard) construction of invariants uses a finite-dimensional
    Lie algebra $\g$ endowed with an invariant non-degenerate scalar product,
     see [BN]. A bit more general construction (see [Ko1]) involves 
     homotopy Lie algebras with scalar products.
     
     In this section we  demonstrate that all these constructions
     are special cases of one universal construction.
     
     The main notion here is the notion of
     a differential $\Z/2\Z$-graded supermanifold,
     or of a $Q$-manifold in short (the terminology is borrowed from [AKSZ],
     the letter $Q$ comes from the standard notation for the generator of
     BRST symmetry in 
     mathematical
     physics).
      By definition, a $Q$-manifold is a super manifold endowed with the action
      of super Lie group $\R^{0|1}$. In other terms, the $Q$-structure
      is given by an odd vector field $Q\in \Pi\Gamma(T_X)$ satisfying
      the equation $[Q,Q]=0$. One defines complex $Q$-manifolds analogously,
       (with possible versions like infinite-dimensional
        manifolds, or partially formal manifolds, like
        our spaces ${\widetilde X}_\C$ as in Section 6).
        
        Basic examples of $Q$-manifolds are
        
        1) $X=\,\,\,$ an ordinary manifold with $Q=0$,
        
        2)$X=\Pi TY=Spec(\bigoplus_i \Omega^i(Y))$, the
         odd tangent space to an ordinary manifold $Y$. Vector field $Q$
         is the de Rham differential,
         
         3) $X=\Pi T_\F Y=Spec(\bigoplus_i \Omega_\F^i(Y))$, an extension of the
          previous example to the case of foliated manifold $(Y,\F)$,
        
        4) $X=\Pi T^{0,1}{\widetilde Y}_{\C}=
        Spec(\bigoplus_i \Omega^{0|i}(Y))$ for complex manifold $Y$, 
        with $Q$ equal to the Dolbeault differential $\db$,
        
        5) $X=\Pi \g=Spec (\bigoplus_i \wedge^i(\g)^*)$, where $g$ is
        a Lie algebra, with the $Q$ equal to the standard differential
         in the cochain complex of $\g$ with trivial coefficients,
         
         6) a tautological extension of the previous example to
          homotopy Lie algebras. Homotopy Lie algebras are defined as
           formal $Q$-manifolds such that the vector field
           $Q$ vanishes at the origin.

    For  $Q$-manifold $X$  
    we define $\Z/2\Z$-graded 
    cohomology group $H^{\bullet}_Q(X)$ 
     as the cohomology $Ker(Q)/Im(Q)$ of the differential
      $Q$ on the super vector space $\O(X)$
     of global functions on $X$ for $C^{\infty}$, Stein , affine,... spaces 
     $X$ (and  sheaf hypercohomology for the non-affine case).
      This cohomology group is equal to
       the space of functions in example 1),
       to the de Rham cohomology in example 2), 
       to the de Rham cohomology for foliations in example 3),
        to the Dolbeault 
      cohomology in example 4),
       and to the Lie algebra cohomology in example 5).

   Instead of flat vector bundles in usual geometry it is convenient
    to speak about $Q$-equivariant vector bundles.
     For example, any flat bundle over a manifold $Y$ produces
      a $Q$-equivariant bundle over $\Pi TY$.
       Any holomorphic bundle over complex manifold $X$ produces
        a $Q$-equivariant bundle over $\Pi T^{0,1}{\widetilde X}_{\C}$.
         Any $\g$-module gives a $Q$-equivariant bundle over $\Pi \g$.
          We define in a uniform way cohomology
           $H^{\bullet}_Q(X,E)$ 
          with coefficients in a $Q$-equivariant bundle $E$ as $Ker(Q)/Im(Q)$
           in the super vector space $\Gamma(X,E)$.
           
           Lie algebras, manifolds, foliations, complex structures, 
           and rational homotopy types
             are all alike.
           
   \vskip 0.5truecm    
 \noindent {9. $Q$-FAMILIES OF SYMPLECTIC MANIFOLDS}    
\vskip 0.5truecm  
     
   The most general construction of characteristic classes including all
    previous cases is the following.
    Let $B$ be a $Q$-manifold and $p:E\ra B$ be  a $Q$-equivariant
     bundle whose fibers are symplectic supermanifolds (may be formal)
      of super dimension $(2n|k)$.
     We assume that the symplectic structure on fibers is also
     $Q$-equivariant. Let $s:B\ra E$ be a $Q$-equivariant section
      of this bundle. The formal completion
      of $E$ along $s(B)$ is a $Q$-equivariant bundle over $B$
       of formal pointed symplectic manifolds.
       Repeating with appropriate modifications constructions from Section 2, 
       we obtain a homomorpism:
       $$H^i_{2n|k}\ra H^{\bullet}_Q(B)\,\,\,.$$

        Here  super vector
   spaces $H^i_{2n|k}$ are defined for any $n,k\ge 0$ starting
   with the standard symplectic super vector space $\R^{2n|k}$ with
    even coordinates $(p_1,\dots,p_n,q_1,\dots,q_n)$ and odd coordinates
      $(\xi_1,\dots,\xi_k)$ and with the symplectic form
      $$\sum_{i=1}^n dp_i\wedge dq_i+\sum_{j=1}^k d\xi_j\wedge d\xi_j\,\,\,.$$

      Graph Cohomology maps to spaces $H^i_{2n|k}$, as in the purely even case.

       The general situation above can be described as a family
        of homotopy Lie algebras with scalar products. Thus, 
        RW invariants is a generalization 
         of the standard construction with homotopy Lie algebras from [Ko1]
          to the case of families.
         
         We repeat here this construction for the case of ordinary 
         even Lie algebras.
       
       If $g$ is a Lie algebra with a non-degenerate invariant scalar product
        then $\Pi \g$ is a flat symplectic super manifold. The vector field 
        $Q$ from example 5) has square equal to $0$ by the Jacobi identity.
         Also, $Q$ preserves the symplectic structure on $\Pi \g$ and vanishes 
         at $0$. The corresponding $Q$-bundle has the base $B$ equal to 
         a point $pt$ with the trivial action of $Q$ on it. The section $s$
          maps $pt$ to $0$.
       
       Analogous constructions can be used in other geometric situations.
        The advantage of symplectomorphism groups is the
         existence of a huge amount of stable classes with possible
          significance for differential topology.

        \vskip 0.5truecm    
 \noindent {10. TOPOLOGICAL QUANTUM FIELD THEORY}    
\vskip 0.5truecm

   Up to now, $3$-dimensional manifolds played  a purely formal role in our 
   exposition. We replaced them from the beginning by cohomology classes
    of symplectomorphism groups. In fact, in [AKSZ] a general Lagrangian
     was constructed, which uses as input data an oriented
     odd-dimensional manifold $M$ 
    and a symplectic manifold $X$. The symmetry group of this Lagrangian 
     is the product of symplectomorphism group 
     $Sympl(X)$ and of certain super extension of the diffeomorphism group
      of $M$. Correlators
       in the corresponding topological
       quantum field theory can be extended to cohomological correlators
        with values in 
        $$H^{\bullet}(BDiff(M),\R)\otimes H^{\bullet}
        (Sympl^{\delta}(X),\R)\,\,\,.$$
        
        It seems that the RW Lagrangian is essentially the same as
        the Lagrangian  in [AKSZ], applied to
      topological quantum field theories depending on parameters. 
      
      In the scheme presented  in [AKSZ],
       one can repace $M$ by an odd-dimensional complex Calabi-Yau manifold.
        The corollary is that finite-type invariants of rational homology
         $3$-spheres 
        give holomorphic
         invariants of $3$-dimensional Calabi-Yau manifolds  with 
         holomorphic volume elements.

       We are planning to discuss it in more details in
        a joint work with A.~Schwarz.

\vskip 0.5truecm

\centerline {\bf Bibliography}\par

\vskip 0.5truecm

\item{[AKSZ]} M.~Alexandrov, M.~Kontsevich, A.~Schwarz, O.~Zaboronsky,
 {\sl The Geometry of the Master Equation and Topological  
 Quantum Field Theory},
  hep-th/9502010.

  \item{[BN]} D.~Bar-Natan, 
  {\sl On the  Vassiliev knot invariants}, Topology, 
  {\bf 34} (1995),  423 - 472.
  
   \item{[BR]} I.~N.~Bernstein,  B.~I.~Rozenfeld,
   {\sl On the characteristic classes of foliations},
    Funkt.\ Anal.\ Appl., {\bf 6} (1972) 1,  68 - 69.
    
    \item{[BH]} R.~Bott, A.~Haefliger,
    {\sl On characteristic classes of $\Gamma$-foliations},
    Bull.\ Am.\ Math.\ Soc., {\bf 78} (1972) No. 6, 1039 - 1044.
   
   \item{[F]} D.~Fuks, {\sl Cohomology of Infinite-Dimensional Lie algebras},
    Consultans Bureau, New York, (1986)
    
    \item{[Ka]} M.~Kapranov, 
    {\sl Rozansky-Witten invariants via Atiyah classes}, alg-geom/9704009.
    
    \item{[Ko1]} M.~Kontsevich, {\sl Feynman diagrams and 
    low-dimensional topology}, First European Congress of 
   Mathematics, 1992, Paris, Volume II, Progress in Mathematics 120,
   Birkhauser 1994, 97 - 121. 
   
 \item{[Ko2]} M.~Kontsevich, {\sl Formal (non)-commutative symplectic
  geometry}, The Gelfand Mathematical Seminars, 1990 - 1992,
   Birkh\"auser (1993), 173 - 187.
   
   \item{[LT]} Le Thang, {\sl An invariant of integral homology
    3-spheres which is universal for all finite type invariants},
     q-alg/9601002.
     
     \item{[RW]} L.~Rozansky, E.~Witten,
      {\sl Hyper-K\"ahler geometry and invariants of $3$-manifolds},
       hep-th/9612216.
   
   \vskip 2truecm
   {\it I.H.E.S., 35 Route de Chartres, Bures-sur-Yvette 91440, FRANCE}
   
   {\it e-mail: maxim@ihes.fr}
\end